\documentclass[a4paper,11pt]{amsart}
\usepackage{hyperref}

\usepackage[usenames,dvipsnames,svgnames,table]{xcolor}
\usepackage[all]{xy}
\usepackage{enumitem}
\usepackage[normalem]{ulem}
\newcommand{\tsk}[1]{\textcolor{YellowOrange}}

\usepackage{tikz}
\usepackage{tikz-cd}
\usepackage{graphicx}
\usepackage{pgf,tikz}
\usetikzlibrary{arrows}
\usepackage[left=2cm,right=2cm, top=3cm, bottom=3cm]{geometry}

\makeatletter
\def\@endtheorem{\endtrivlist}
\makeatother

\usepackage[active]{srcltx}
\usepackage{amsmath}
\usepackage{amssymb}
\usepackage{amscd}
\usepackage{amsthm}
\usepackage{mathrsfs}  
\usepackage{mathtools, nccmath}

\usepackage{nicefrac}

\newcommand{\Pic}{\operatorname{Pic}}
\newcommand{\Nef}{\operatorname{Nef}}

\newtheorem{teo}{Theorem}[section]

\newtheorem*{mainteo}{Main Theorem}

\newtheorem{defin}[teo]{Definition}
\newtheorem{proposition}[teo]{Proposition}
\newtheorem{cor}[teo]{Corollary}

\newtheorem{lemma}[teo]{Lemma}

\theoremstyle{definition}

\newtheorem{remark}[teo]{Remark}

\newtheoremstyle{dico}
{\baselineskip}   
{\topsep}   
{}  
{0pt}       
{} 
{.}         
{5pt plus 1pt minus 1pt} 
{}          
\theoremstyle{dico}

\numberwithin{equation}{section}

\newcounter{example}[subsection]

\newcommand{\ra}{\rightarrow}

\newcommand{\Zeta}{{\mathbb{Z}}}

\newcommand{\restr}[1]          {\vert_{#1}}

\newcommand{\om}{\omega}

\renewcommand{\phi}{\varphi}

\newcommand{\OO}{\mathcal{O}}

\newcommand{\Hilb}{S^{[2]}}
\newcommand{\projo}{\mathbb{P}(\Omega^1_S)}

\newcommand{\Sym}{{\operatorname{Sym}}}

\newcommand{\mihi}[1]{}

\newcommand{\e}{\qquad \text{and}\qquad}

\begin{document}

\pagestyle{myheadings}

\title{Gauss-Prym maps on Enriques surfaces}

\author{Dario Faro and Irene Spelta}

\address{Universit\`a degli Studi di Pavia, Dip. di Matematica, Via Ferrata 5, 27100 Pavia, Italy.}
\address{Dario Faro }
\email{d.faro1@campus.unimib.it}

\address{Irene Spelta }
\email{irene.spelta@unipv.it}

\begin{abstract}{We prove that the $k$-th Gaussian map $\gamma^k_{H}$ is surjective on a polarized unnodal Enriques surface $(S, H)$ with $\phi(H)>2k+4$. In particular, as a consequence, when $\phi(H)>4(k+2)$, we obtain the surjectivity of the $k$-th Gauss-Prym map $\gamma^k_{\om_C\otimes\alpha}$, with $\alpha:=\omega_{S\restr{C}}$, on smooth hyperplane sections $C\in \vert H\vert.$ In case $k=1$ it is sufficient to ask $\phi(H)>6$.}
\end{abstract}

\maketitle

\section{Introduction}
In this paper, we show the surjectivity of the $k$-th Gauss-Prym map $\gamma^k_{\om_C\otimes\alpha}$, with $\alpha:=\omega_{S\restr{C}}$, on smooth hyperplane sections of polarized unnodal Enriques surfaces.

Let us briefly introduce the setting. Let $X$ be a smooth projective variety and $L,M$ line bundles on it. We denote by $\mu_{L,M}$ the multiplication map on global sections and by $K(L,M)$ its kernel. The  first Gaussian map associated to $L$ and $M$ is the map \[\gamma^1_{L,M}: K(L,M)\ra H^0(X, \Omega^1_X\otimes L\otimes M). \]
locally defined as $\gamma^1_{L,M}(s\otimes t)=sdt-tds$. As usual, set $\gamma^1_{L}:=\gamma^1_{L,L}$. Moreover, since $\gamma^1_{L}$ vanishes identically on symmetric tensors, we usually write $\gamma^1_{L}:\Lambda^2H^0(L)\ra H^0( \Omega^1_X\otimes L^2)$. 

The first Gaussian map  $\gamma^1_{L,M}$ can be seen as the first brick of a hierarchy of maps, the higher Gaussian maps $\gamma^k_{L,M}$:
\[\gamma^k_{L,M}: H^0(X\times X, \mathcal{I}_{\Delta_X}^k(L\boxtimes M) )\ra H^0(X,\Sym^k\Omega^1_X\otimes L\otimes M ),\]
where $\mathcal{I}_{\Delta_X}$ is the ideal of the diagonal $\Delta_X\subset X\times X$ (\cite{wahl jac},  for more details see section \ref{sezione 1}).

Gaussian maps turn out to be useful to study the geometry of $X$ once one has a good understanding of the behaviour of its hyperplane sections. In particular, they can be used to distinguish those curves that are hyperplane sections of $X$. 

The prototype for this kind of argument goes back to Wahl. In \cite{wahl jac} he proved that if $C$ is a hyperplane section of a K3 surface, then the Gaussian map $\gamma^1_{\om_C}$ (sometimes referred to as the Wahl map) is not surjective. 
The significance of this result becomes more evident once compared with the result of Ciliberto, Harris, and Miranda (\cite{chm}) stating that the Wahl map of the generic genus $g$ curve is surjective as soon as this is numerically possible, i.e. for $g\geq 10$ with the exception of $g=11$. On the other hand for $g<10$ and $g=11$ it is known that the generic curve lies on a K3 surface (\cite{mm}). 

It is worth to mention that the result of Wahl on $\gamma^1_{\om_C}$ has been independently obtained also by Beauville and Merindol (\cite{bm}). While \cite{wahl jac} relates the corank of $\gamma^1_{\om_C}$ with deformations of the cone over $C$, \cite{bm} links the same corank with the splitting of the normal bundle exact sequence of $C$ in $X$.

Note that quite recently, Arbarello, Bruno and Sernesi gave in \cite{abs} a partial converse to  \cite[Theorem 5.9]{wahl jac}, \cite[Proposition 5]{bm}. Indeed, they showed that Brill-Noether-Petri general curves lying as hyperplane sections on K3 surfaces (or on singular specialization of them) are exactly those curves where the Wahl map $\gamma^1_{\om_C}$ is not surjective. 

In a similar fashion, Colombo, Frediani and Pareschi investigated abelian surfaces in \cite{cfp}. They proved, using a ``Beauville-Merindol" approach, that the first Wahl map on hyperplane sections is ``tendentially" surjective, whereas the second one is never. 

Also the case of Enriques surfaces has been treated by some authors. In  \cite[Corollary 1.2]{cv}, Ciliberto and Verra  proved that  $\gamma^1_{\om_C, \om_C\otimes \alpha}$ is not surjective if $C$ a hyperplane section of an Enriques surface and  $\alpha:=\om_{S\restr{C}}$. 
Again, this becomes more interesting once compared with \cite[Theorem 1.5]{cv} which proves $\gamma^1_{\om_C, \om_C\otimes \alpha}$ to be generically surjective for $g\geq 12, g\neq 13$. Notice that, in this case, the generality is required in $\mathcal{R}_g$, namely in the (coarse) moduli space of pairs $[C,\alpha]$, where $C$ is a smooth projective genus $g$ curve and $\alpha$ is a non-trivial 2-torsion line bundle on it. It is well-known that the pair $[C,\alpha]$ identifies a principally polarized abelian variety, the Prym variety. This is the reason why the line bundle $\om_C\otimes \alpha$ is usually called Prym-canonical line bundle.

We finally mention that the study of  the surjectivity of Gaussian maps of the form $\gamma^1_{M, \omega_C}$ for $C$ curves on Enriques surfaces and  $ M$ line bundle on $C$, carried out by Knutsen and Lopez in \cite{knudue}, has been applied  in \cite{knutre} to study threefolds whose general hyperplane section is an Enriques surface. 

The goal of this paper has the same flavour of the aforementioned papers. Indeed, we investigate smooth hyperplane sections $C$ of projective Enriques surfaces $S\subseteq\mathbb{P}^g$ with respect to the $k$-th Gaussian (Gauss-Prym) map: \begin{equation}\label{mappa gauss sotto}
	\gamma^k_{\om_C\otimes \alpha}: H^0(C\times C, \mathcal{I}^k_{\Delta_C}((\om_C\otimes \alpha)\boxtimes (\om_C\otimes \alpha))\ra H^0(C, \om_C^{\otimes k+2}).
\end{equation} 
To this purpose, we keep in mind that Barchielli and Frediani in \cite{bf}, respectively Colombo and Frediani in \cite{cf3}, proved that the surjectivity holds for the first, respectively the second, Gauss-Prym for the general point $[C,\alpha]\in \mathcal{R}_g$ provided that $g\geq 12$, respectively $g\geq 20$.

In order to investigate the Gauss-Prym map \eqref{mappa gauss sotto}, we put it in the diagram \eqref{diagramma} and we study its surjectivity showing that it holds for the other three involved maps. 

The first step focuses on the $k$-th Gaussian map (on $S$) \begin{equation}\label{gaussiana sopra}
	\gamma^k_{\OO_S(C)}: H^0(S\times S, \mathcal{I}^k_{\Delta_S}(\mathcal{O}_S(C)\boxtimes\mathcal{O}_S(C)))\ra H^0(S, \Sym^k\Omega^1_S\otimes \OO_S(2C)). 
\end{equation}
Here, a very recent paper by Rios Ortiz (\cite{ortiz}) is fundamental. There the author studies higher Gaussian maps $\gamma^k_{L}$. In particular, motivated by previous works on $\gamma^2_{\om_C}$ (\cite{ccm,cf1,cf2}), he shows that the maps $\gamma^k_{\om_C}$ ($ k>1 $) are surjective on hyperplane sections of K3 surfaces. His proof relates the surjectivity  with the cohomology of certain linear systems on the Hilbert scheme of two points on  a projective and regular surface  $S$. More precisely, he shows that if $H^1(S^{[2]}, L-(k+2)B)=0$ then $\gamma^k_{L}$ is surjective. Here $2B$ is the exceptional divisor of the Hibert-Chow morphism $S^{[2]}\ra S^{(2)}$ that resolves the singularities of the symmetric product $S^{(2)}$. Notice that the descriptions both of the Picard group of $S^{[2]}$, due to Fogarty (\cite{fog1,fog2}), and of the cone of its nef divisors, due to Nuer (\cite{nuer1}) in case of an unnodal Enriques surface $S$,  are crucial for our analysis. We will recall them accordingly. It turns out that much of the geometry of $S$, and hence of $S^{[2]}$, is governed by a function, the $\varphi$-function, defined on $\Pic(S)$ and taking values in $\Zeta_{\geq0}$ (see Definition \ref{phi}). Our first result is the following:
\begin{teo}\label{teo 1}
	Let $S$ be an unnodal Enriques surface, $H$ be a line bundle on $S$ with  $\varphi(H)>  2k+4$ and $C\in \vert H\vert$. The $k$-th Gaussian map $\gamma^k_{\OO_S(C)}$ is surjective. 
\end{teo}
%

As second step, we study the restriction map \[p_1: H^0(S, \Sym^k\Omega^1_S\otimes \OO_S(2C))\ra H^0(C,\Sym^k \Omega^1_S\otimes \OO_S(2C)\restr{C}).\]
The surjectivity is proved identifying the exceptional divisor $2B$ of the Hilbert-Chow morphism with $\mathbb{P}(\Omega_S^1)$ and using properties of projective bundles. Moreover, we show that the case $k=1$ can be solved more easily using \cite[Lemma 6.3]{knutsenetal}. Indeed, in this (very exhaustive) paper the authors systematically investigate curves on Enriques surfaces. 


Finally, simply from the $k$-th symmetric power of the normal bundle sequence of $C$ in $S$, we show the surjectivity of \[p_2:H^0(C,\Sym^k \Omega^1_S\otimes \OO_S(2C)\restr{C})\ra H^0(C, \om_C^{\otimes k+2}). \] 

In conclusion, this procedure leads to our 
\begin{mainteo}
	Let $C$ be a smooth hyperplane section of an unnodal Enriques surface $(S,H)$ with $\phi(H)>4(k+2)$. Then the $k$-th Gauss-Prym map $\gamma^k_{\om_C\otimes \alpha}$ is surjective. In case $k=1$ it is sufficient to ask $\phi(H)>6$.
\end{mainteo}
Notice that, in view of the results of \cite{bf, cf3}, the first two Gauss-Prym maps do not give a distinguishing criterion for curves on (unnodal) Enriques surfaces. \\

The paper is organised as follows. In section \ref{sezione 1}, we review some preliminaries on Gaussian maps and the recent result of Rios Ortiz. In section \ref{sezione 2}, we focus on Enriques surfaces. In particular, we study positivity properties of line bundles on Enriques surfaces. Then we recall Nuer's description of $\Nef(\Hilb)$. Finally, in section \ref{sezione 3}, we prove our Main Theorem.  \\

\section{Gaussian maps and Hilbert Scheme of points}\label{sezione 1}
In this section we introduce Gaussian maps and we explain how to reinterpret them in terms of the cohomology of certain line bundles on the Hilbert scheme of points. In particular, we refer to results of \cite{wahl jac,wahl intro} and of \cite{ortiz}.

Let $X$ be a smooth projective variety, $\Delta_X\subset X\times X$ the diagonal, and $\mathcal{I}_{\Delta_X}$ its ideal sheaf. Moreover, let $L,M$ be two line bundles on $X$. On $X\times X$ we consider the external product $L\boxtimes M:=\pi_1^*L\otimes \pi_2^*M$, where $\pi_i$ is the projection to the $i$-th factor. For any non-negative integer $k$, there is a short exact sequence \begin{equation}\label{ideale}
	0\ra \mathcal{I}^{k+1}_{\Delta_X}\ra \mathcal{I}^k_{\Delta_X}\ra \Sym^k\Omega^1_X\ra 0
\end{equation} 

\begin{defin}
	The $k$-th Gaussian map $\gamma^k_{L,M}$ is the map on global sections \[\gamma^k_{L,M}: H^0(X\times X,\mathcal{I}^k_{\Delta_X}(L\boxtimes M) )\ra H^0(X,\Sym^k\Omega^1_X\otimes L\otimes M )\] induced by the exact sequence \eqref{ideale} twisted by $L\boxtimes M$. 
\end{defin}
Notice that, by definition, for each $k\geq 1$, the domain of $\gamma^k_{L,M}$ is the kernel of the previous one:
\[\gamma^k_{L,M}: \ker(\gamma^{k-1}_{L,M})\ra H^0(X,\Sym^k\Omega^1_X\otimes L\otimes M ). \]
We will assume that the two line bundles $L$ and $M$ coincide. In this case the standard notation is $\gamma^k_{L,L}=:\gamma^k_{L}$.  Here, in particular, we will deal with the surjectivity of the $k$-th Gaussian map $\gamma^k_L$ defined on an Enriques surface $S$ with line bundles $L=\OO_S(H)$, under suitable assumptions on $H$.


Let us assume $S$ to be a smooth projective surface. Now we discuss how the Hilbert scheme of two points on $S$ can be used to study the surjectivity of the $k$th Gaussian map. We borrow from \cite{ortiz} what follows. 

Let $S^{(2)}$ be the second symmetric product of $S$ and  $S^{[2]}$ be the Hilbert scheme of $2$ points on $S$. By Fogarty (\cite{fog1,fog2}), $S^{[2]}$ is smooth and isomorphic to the blow-up of $S^{(2)}$ along the diagonal $\Delta_S$ of $S^{(2)}$. Let us denote by $2B$ the smooth exceptional divisor. Moreover, if $H^1(S,\mathcal O_S)=0$, then \cite[Theorem 6.2]{fog2} shows
\begin{equation}\label{gruppo picard Fogarty}
	\Pic(S^{[2]}) \simeq \Pic(S) \oplus \mathbb{Z}B.
\end{equation}
Indeed, let 
$$
\rho: S^{[2]} \rightarrow S^{(2)}
$$
be the  Hilbert-Chow morphism (i.e. the blow-up morphism). Let $\pi_{i}: S \times S \ra S, \ i=1,2$ be the two projections. If $L \in \Pic(S)$, then $L\boxtimes L=\pi_{1}^*L \otimes \pi_{2}^*L$ defines a $\mathbb Z/2\mathbb Z$-invariant line bundle which naturally descends to a line bundle over $S^{(2)}$. Pulling it back through $\rho$, we obtain a line bundle  $\tilde L$ over $S^{[2]}$. This procedure gives an injective group homomorphism of  $\Pic(S)$ into  $\Pic(S^{[2]})$. In particular, we will make use of the fact that
\begin{equation}\label{canonico}
\widetilde{K_S}=K_{S^{[2]}}.
\end{equation}
Indeed, see e.g. \cite{gotsche}, the canonical divisor $K_{S^2}$ is invariant under the action of $\mathbb{Z}/2\mathbb{Z}$. Thus it gives a canonical Cartier divisor $K_{S^{(2)}}$ on $S^{(2)}$. Being $\rho$ crepant (\cite{beau}), we have $\rho^*(K_{S^{(2)}})=K_{S^{[2]}}$. Thus, via the identification \eqref{gruppo picard Fogarty}, we get $\eqref{canonico}$.

In \cite{ortiz}, Rios Ortiz shows the following:
\begin{proposition}[\cite{ortiz}, Corollary 2.7]
	For every $p, k\geq0$ there are natural identifications\[H^p(S\times S, \mathcal{I}^k_{\Delta_S}(L\boxtimes L))\cong H^p(\Hilb, \tilde L -kB)\oplus  H^p(\Hilb, \tilde L -(k+1)B). \]
\end{proposition}
This is crucial in providing the following criterion for the surjectivity of higher Gaussian maps. Indeed, in loc. cit., Rios Ortiz states the following:
\begin{teo}[\cite{ortiz}, Theorem 2.9]\label{ortiz}
	If $H^1(\Hilb, \tilde L -(k+2)B)=0$, then $\gamma^k_L$ is surjective. 
\end{teo}
\section{Line bundles on Enriques surfaces}\label{sezione 2}

An Enriques surface is  a smooth complex projective surface $S$ such that $h^1(S, \mathcal O_S)=0$,  $\om_S^2=\mathcal O_S$ but $\om_S \neq 0$. An Enriques surface is said to be unnodal if it does not contain $-2$ curves, that is, irreducible curves with self intersection $-2$ (necessarily isomorphic to $\mathbb{P}^1$). If $H$ is an ample line bundle on an Enriques surface such that $H^2 =2g-2$ then  dim$(|H|)=g-1$ and the general element of $|H|$ is a smooth irreducible curve of genus $g$. 

The geometry of curves on Enriques surfaces is strongly related to the so called $\phi$-function introduced by Cossec.
\begin{defin}
	\label{phi}
	Let $H$ be a line bundle on an Enriques surface $S$ such that $H^2 >0$.  
	$$
	\phi(H):=\text{min}\{|H \cdot F|: F\in \Pic(S), F^2=0, F \not \equiv 0\}.
	$$
\end{defin}
In what follows, we need the relation between $\phi(H)$ and the ``regularity" of the map to the projective space induced by $H$. For this purpose, we recall the notion of $k$-very ampleness.
\begin{defin}
	Let $S$ be a smooth connected surface over the complex numbers and let $k \geq 0$ be an integer. A line bundle $H$ is said to be $k$-very ample if for any $0$-dimensional subscheme $(Z, \mathcal O_Z)$ of length $k+1$ the restriction map $H^0(H) \rightarrow H^0(H \otimes \mathcal O_Z)$ is surjective.   
\end{defin}
\begin{remark}
	Observe that a line bundle is $0$-very ample if and only if it is globally generated and $1$-very ample if it is very ample. 
\end{remark}
Knutsen and Szemberg proved independently the following:
\begin{teo}[\cite{szem}, \cite{knuno}]
	Let $S$ be an  Enriques surface. Then $H$ is $k$-very ample if and  only if $\phi(H) \geq k+2$ and there exists no effective divisor $E$ such that $E^2=-2$ and $H \cdot E \leq k+1$.  
\end{teo}
	%
	%
	%
	

	Thus, one immediately gets the following
	\begin{cor}\label{k very ampiezza}
		Let $S$ be an unnodal Enriques surfaces. The line bundle $H$ is $k$-very ample if and only if $\phi(H) \geq k+2$.
	\end{cor}
	
	The notion of $k$-very ampleness is useful to construct (very) ample line bundles on Hilbert scheme of points. Indeed, let $H$ be a line bundle on a smooth projective surface $S$ and let $Z$ be a $0$-dimensional subscheme of $S$ of length $k+1$. Let
	$$
	0 \rightarrow H \otimes \mathcal{I}_Z \rightarrow H \rightarrow H \otimes \mathcal O_Z \rightarrow 0.
	$$
	be the exact sequence defining $Z$ as a subscheme, tensored by $H$. 
	If $H$ is $k$-very ample, $H^0(S,H \otimes \mathcal{I}_Z)$ is a codimension $k+1$ linear subspace of  $H^0(S, H)$. Thus, we have a map:
	\begin{align}\label{mappa nella grass}
		\psi_H: S^{[k+1]} &\rightarrow Gr((k+1),h^0(S,H)) 
		\\ 
		(Z,\mathcal{I}_Z) &\rightarrow H^0(S,H )/ H^0(S,H \otimes \mathcal{I}_Z).\notag
	\end{align}
	In \cite{catanesegot}, Catanese and G\"oettsche showed that \eqref{mappa nella grass} is an embedding if and only if $H$ is $(k+1)$-very ample. Moreover, it is well-known that $\tilde H-B$ is the pull-back of the very ample line bundle $\mathcal O(1)$ on the Grassmanian, see e.g. \cite{bertram coskun}. Therefore, if $H$ is $(k+1)$-very ample then $\tilde H-B$ is very ample.  For our convenience, we state the following 
	\begin{proposition}\label{h-b very ample}
		Let $S$ be an unnodal Enriques surface and $H$ a line bundle such that $\phi(H) \geq 4$. Then $\tilde H - B$ is a very ample line bundle on $\Hilb$.
		\begin{proof}
			From Corollary \ref{k very ampiezza}, we have that $H$ is $2$-very ample, Then the statement follows since $\psi_H$ is an embedding.
		\end{proof}
	\end{proposition}
	
	We end this section recalling the description of the cone of nef divisors $\Nef(\Hilb)$  given by Nuer in \cite{nuer1}.
	\begin{teo}[\cite{nuer1}]\label{teo nuer}
		Let $S$ be an unnodal Enriques surface and $k\geq 2$. Then $\tilde L -aB\in \Nef(S^{[k]})$ if and only if $L\in \Nef(S)$ and $0\leq a\leq \phi(L)/k.$ 
		
	\end{teo}

	\section{Proof of the Main Theorem}\label{sezione 3}
	Let $S$ be an unnodal Enriques surface,   $k \geq 1$. Let  $H\in \Pic(S)$ with $\phi(H)>4(k+2)$ and $C\in \vert H\vert$ be a smooth hyperplane section. Then $\alpha:=\om_{S\restr{C}}$ defines a 2-torsion line bundle on $C$ and so $\om_C\otimes \alpha$ is Prym-canonical on $C$. 
	
	In this section we will show that the $k$-th Gauss-Prym map  
	\begin{equation} \label{gaussiana sotto}
		\gamma^k_{\om_C\otimes \alpha}:H^0(C\times C,\mathcal{I}^k_{\Delta_C}((\omega_C \otimes \alpha) \boxtimes (\omega_C \otimes \alpha) ) )\ra H^0(C,\omega_C^{\otimes k+2})
	\end{equation}
	is surjective. We look at the following diagram 
	\begin{equation}\label{diagramma}
		\begin{tikzcd}
			H^0(S\times S,\mathcal{I}^k_{\Delta_S}( \mathcal O_S(C)\boxtimes \mathcal O_S(C) ) )	\arrow{dd}\arrow{r}{\gamma^k_{\mathcal{O}_S(C)}}&  H^0(S, \Sym^k  \Omega^1_S(2C))\arrow{dr}{p_1}&\\
			&  & H^0(C, \Sym^k  \Omega^1_S(2C)\restr{C})\arrow{dl}{p_2}\\
			H^0(C\times C,\mathcal{I}^k_{\Delta_C}((\omega_C \otimes \alpha) \boxtimes (\omega_C \otimes \alpha )) )\arrow{r}{\gamma^k_{\om_C\otimes \alpha}} &H^0(C,\omega_C^{\otimes k+2}).&
		\end{tikzcd}
	\end{equation}
	Here $\gamma^k_{\mathcal{O}_S(C)}$ is the $k$-th Gaussian map (on $S$) associated with the line bundle $\mathcal{O}_{S}(C)$. The vertical arrow and $p_1$ are restriction maps. Finally, $p_2$ comes from the $k$-th symmetric power of the conormal bundle sequence
	\begin{equation}\label{symmetric conormal}
		0\ra \Sym^{k-1} \Omega^1_{S_{\restr{C}}}(-C)\ra \Sym^{k}\Omega^1_{S_{\restr{C}}} \ra\omega_C^{\otimes k} \ra 0
	\end{equation}
	tensored by $\OO_C(2C)$.
	
	We prove that $\gamma^k_{\mathcal{O}_S(C)}, p_1, \text{and } p_2$ are surjective. From this we obtain the surjectivity of \eqref{gaussiana sotto}. \\
	
	{\bfseries \noindent{Surjectivity of $\gamma^k_{\mathcal{O}_S(C)}$}}. Take $S$, $H \in \Pic(S)$ such that $\phi(H)>2(k+2)$ and $C$ as above (notice that the assumption on $\phi(H)$ is weaker for this step).  By Theorem \ref{ortiz}, to show the surjectivity of $\gamma^k_{\OO_S(C)}$, we just need to show that 	$H^1(S^{[2]}, \tilde H -(k+2)B)=0$ and to show the latter it is by \eqref{canonico} and Kodaira vanishing enough to prove
	that $ \widetilde{H-K_S}-(k + 2)B $ is ample. Let us start with the following:
	\begin{proposition}
	    
	\end{proposition}\label{line bundles ampi}
		Let $S$ be an unnodal Enriques surface. If $L$ is a line bundle on $S$ such that $\phi(L)> 2k+4$, then $\widetilde{L}-(k+2)B$ is ample on $\Hilb$.
		\begin{proof}
			Since S is unnodal, L is nef. Therefore, using Nuer's description of $\Nef(\Hilb)$ (see Theorem \ref{teo nuer}), we conclude that 
			$$\widetilde{L}-(k+1)B\e \widetilde{L}-\frac{\phi(L)}{2}B$$
			both belong to $\Nef(\Hilb)$. Furthermore, for the same reasons,  $\widetilde{L}-(k+2)B$ is nef. 
			
			Assume now, by contradiction, that there exists $D\in \overline{NE}(S^{[2]})$ violating Kleiman's criterion for ampleness, namely such that $(\widetilde{L}-(k+2)B)\cdot D=0$. Thus we have $\widetilde{L}\cdot D=(k+2)(B\cdot D)$. 
			Since \[(\widetilde{L}-(k+1)B)\cdot D\geq 0 \e (\widetilde{L}-\frac{\phi(L)}{2}B)\cdot D\geq 0,   \]
			we obtain \[B\cdot D=\tilde L\cdot D=0.\]
			Indeed, $(\widetilde{L}-(k+1)B)\cdot D\geq 0$ yields $(k+2)(B\cdot D)=\widetilde{L}\cdot D\geq (k+1)(B\cdot D)$, thus $B\cdot D\geq 0$. On the other hand, $ (\widetilde{L}-\frac{\phi(L)}{2}B)\cdot D\geq 0 $ yields $(\frac{\phi(L)}{2}-k-2)(B\cdot D)\leq 0$, thus $B\cdot D\leq 0.$
			
			By Proposition \ref{h-b very ample}, $\tilde L-B$ is very ample. Thus, the condition \[(\tilde L-B)\cdot D>0\; \text{for all} \; D\in \overline{NE}(S^{[2]}) \] yields a contradiction.

		\end{proof}
	\begin{cor}
		\label{coroampleds}
		Let $H\in \Pic(S)$ be such that $\phi(H)>2k+4$. Then $\widetilde{H-K_S}-(k+2)B$ is ample. 
		\begin{proof}
			Being $K_S$ numerically trivial, we have $\phi(H-K_S)=\phi(H)$. Therefore, we just need to apply Proposition \ref{line bundles ampi}.
		\end{proof}
	\end{cor}
	This ends the proof of Theorem \ref{teo 1}.\\

	\begin{remark}
		Notice that the strategy of Proposition \ref{line bundles ampi} can be used to produce some ample line bundles on $\Hilb$: let $L$ be a line bundle on $S$ such that $\phi(L)=k, k>4$. We claim that \begin{equation}\label{altri divisori ampi}
			\tilde L-\bigg(\frac{k}{2}-1-r\bigg)B \quad \text{is ample for }\; 1\leq r<\frac{k}{2}-1.
		\end{equation} 
		Indeed, by Theorem \ref{teo nuer}, we have $\tilde L-\frac{k}{2}B$ and $\tilde L-(\frac{k}{2}-2)B$ in $\Nef(\Hilb)$. Arguing as before by contradiction, we obtain $\tilde L-(\frac{k}{2}-1)B$ ample. This implies \eqref{altri divisori ampi}. Indeed, if there exists $D\in \overline{NE}(S^{[2]})$ such that $(\tilde L-(\frac{k}{2}-1-r)B)\cdot D=0$, then $$\bigg(\tilde L-\bigg(\frac{k}{2}-1\bigg)B\bigg)\cdot D=-r(B\cdot D)$$ and so, being the left hand side strictly positive, we would get $B\cdot D<0$. Now, since $\tilde L$ is nef, we would have \[0\leq\tilde L\cdot D= \bigg(\frac{k}{2}-1-r\bigg)(B\cdot D)<0\]
		when $1\leq r<\frac{k}{2}-1$. Thus we have a contradiction. \\
		
	\end{remark}

	{\bfseries \noindent{Surjectivity of $p_1$}}. Let $S$, $H \in \Pic(S)$ be such that $\phi(H)>4(k+2)$ and $C \in |H|$. Let us consider the following short exact sequence:
	\[0\ra \Sym^k \Omega^1_S(C) \ra \Sym^k\Omega_S^1(2C)\ra \Sym^k \Omega^1_{S\restr{C}}(2C)\ra 0. \]
	In order to prove the surjectivity of $p_1$, it is enough to prove the following lemma.
	\begin{lemma}
	\label{piunocondiann}
	Let $S$ be an unnodal Enriques surface and $H \in \Pic(S)$ such that $\phi(H) >4(k+2)$. Then
	\begin{equation*}
		H^1(S,\Sym^k \Omega^1_S(C))=0 \; \ \text{for all} \; k \geq0.
	\end{equation*}
	\end{lemma}
	\begin{proof}
	The case $k=0$ follows from the exact sequence
	\[0\ra \omega_S(-C) \ra \omega_S\ra \omega_{S\restr{C}}\ra 0. \]
	 using Serre duality and the fact that $S$ is an Enriques surface. For $k \geq 1$ we proceed in the same way as in \cite{ortiz}. Let $
	\pi: \mathbb{P}(\Omega^1_S) \ra S$ be the projectivisation of  $\Omega^1_S$ and let $\xi$ be the class of the tautological line bundle $ \mathcal O_{\mathbb{P}(\Omega^1_S)}(1)$ on $\mathbb{P}(\Omega^1_S)$. Then (see for example  \cite[Appendix A]{laz2}) the following properties hold:
	\begin{align}
		\label{lazarproj}
		\pi_{*}(\mathcal O_{\projo} (k\xi))&=\Sym^k \Omega^1_S; 
		\\
		R^i \pi_{*} (\mathcal O_{\projo} (k\xi))&=0 \ \  \forall i >0. \nonumber
	\end{align}
	By the projection formula, it then follows that 
	$$
	R^i \pi_{*} (\mathcal O_{\projo} (k\xi + \pi^*H))=R^i \pi_{*} (\mathcal O_{\projo} (k\xi)) \otimes \mathcal O_S(C)=0 \ \  \forall i >0,
	$$
	and so, by degeneration of the Leray spectral sequence, one gets
	$$
	H^1(\mathbb{P}(\Omega^1_S),\mathcal O_{\projo} (k\xi  +\pi^*H)) \simeq H^1(S,\pi_*(\mathcal O_{\projo}(k\xi + \pi^*H))) \simeq H^1(S,\Sym^k \Omega^1_S(C)).
	$$
	
	The bundle $\mathbb{P}(\Omega^1_S)$ embeds in $S^{[2]}$ as the exceptional divisor $2B$. Under this identification, one has 
	\begin{equation}
		\label{equkawasurj}
		\tilde H_{\restr{2B}}=2\pi^*H \e 
		\xi=-B_{\restr{2B}}. 
	\end{equation}
	To show that $H^1(\mathbb{P}(\Omega^1_S),k \xi +\pi^*H)=0$, we apply the Kodaira vanishing theorem. We prove that $ k \xi + \pi^*H - K_{\mathbb{P} (\Omega^1_S)}$ is ample. Using \eqref {equkawasurj}, \eqref{canonico}, and the adjunction formula for the divisor $\mathbb{P}(\Omega^1_S) \simeq 2B$, we get 
	$$
	2(k \xi + \pi^*H - K_{\mathbb{P}(\Omega^1_S)})=(-2kB + \widetilde{H}-2\tilde K_S - 4B)_{|_{2B}}=(\widetilde{H} - 2(k+2)B)_{|_{2B}}.
	$$
	The assumption $\phi(H) > 4(k+2)$ allows to conclude. Indeed, by Proposition \ref{line bundles ampi}, the latter is the restriction of an ample line bundle to a smooth divisor. Hence it is ample. 
	\end{proof}

	

		{\bfseries \noindent{Surjectivity of $p_2$.}} Assume $H \in \Pic(S)$ with $\phi(H)>4(k+1)$. By \eqref{symmetric conormal} twisted by $\OO_C(2C)$, it is enough to show that 
		$ H^1(C,\Sym^{k-1} \Omega^1_{S_{|_{C}}}(C))=0, $  equivalently \begin{equation}\label{condizione p_2}
		H^0(C,\Sym^{k-1} \mathcal{T}_{S_{|_{C}}}(\alpha))=0,
		\end{equation}  by Serre duality. By
		\[0\ra \Sym^{k-1} \mathcal{T}_{S}(-C+ K_S)\ra \Sym^{k-1} \mathcal{T}_{S}(K_S)\ra \Sym^{k-1} \mathcal{T}_{S_{|_{C}}}(\alpha) \ra 0,\]
	    it is enough to prove that
	\begin{equation}\label{annullamento per p_2}
H^0(S,\Sym^{k-1} \mathcal{T}_{S}(K_S))=0 \e H^1(S,\Sym^{k-1} \mathcal{T}_{S}(-C+K_S))=0.
	\end{equation}
The right hand vanishing follows from Serre duality and Lemma \ref{piunocondiann}, along with our assumptions on $\phi(H)$. To prove the left hand vanishing of \eqref{annullamento per p_2}, let $Y \xrightarrow{\pi} S$ be the $K3$ double cover of $S$, namely the degree 2 (cyclic) covering associated with the pair $(S, \om_S)$. As $\pi$ is unramified, we  have $\Sym^{k-1} \mathcal{T}_{Y} =\pi^*\Sym^{k-1} \mathcal{T}_{S}$. Using the  projection formula we obtain
		\begin{align}
			\nonumber
			H^0(Y,\Sym^{k-1} \mathcal{T}_{Y})
			=H^0(S,\Sym^{k-1} \mathcal{T}_{S}) \oplus H^0(S,\Sym^{k-1} \mathcal{T}_{S}(K_S)) \nonumber.
		\end{align}
		A result of Kobayashi (\cite{kob}) asserts that $H^0(Y,\Sym^{k-1} \mathcal{T}_{Y})=0$, and so we get $H^0(\Sym^{k-1} \mathcal{T}_{S}(K_S))=0$, as desired. 
		This ends the proof of the Main Theorem. 
		
		\begin{remark}
			\label{punocasok1}
			The case $k=1$ can be treated in a easier way. Notice that it is sufficient to require $\phi(H)> 6$ (instead of $\phi(H)> 12$). Indeed, by Theorem \ref{teo 1}, if $\phi(H)> 6$, the map $\gamma^1_{\mathcal O_S(C)} $ is surjective.  
			As for $p_1$, when $k=1$, condition \eqref{piunocondiann} becomes $$
			H^1(S, \Omega^1_S(C))=0,
			$$
			which by Serre duality is equivalent to \begin{equation}\label{coso che si deve annullare}
				H^1(S, \mathcal T_S\otimes\om_S(-C))=0.
			\end{equation}
			Let $Y \xrightarrow{\pi} S$ be the $K3$ double cover of $S$ as above.  Again, using the projection formula, we get   
			\begin{equation}\label{decomposition p_1}
				H^1(Y,\mathcal T_{Y}(-(\pi^*H))=H^1(Y,\pi^*(\mathcal T_{S}(-H)))=H^1(S,\mathcal T_{S}(-H)) \oplus H^1(S,\mathcal T_{S}\otimes\om_S(-H)).
			\end{equation}			
			By \cite[Lemma 6.3.]{knutsenetal}, if $\phi(H)\geq5$ then  
			$
			H^1(Y,\mathcal T_{Y}(-\pi^*H))=0,
			$
			thus, using \eqref{decomposition p_1}, we get \eqref{coso che si deve annullare}. 
			This yields the surjectivity of $p_1$. 
			
			The surjectivity of $p_2$ follows as \eqref{condizione p_2} for $k=1$ reads $H^0(C,\alpha)=0$, which is satisfied since $\alpha$ is a non-trivial 2-torsion line bundle on $C$.
		\end{remark} 
	
\section*{Acknowledgments}
We kindly thank Paola Frediani for introducing us to the topic of this paper and for the precious conversations we had. Moreover, the first author would like to thank Angel D. Rios Ortiz for a valuable discussion. We also thank Thomas Dedieu for catching an inaccuracy in our first draft. Finally, we thank the two referees for their very useful comments.

	The authors were partially supported by MIUR PRIN 2017
``Moduli spaces and Lie Theory'',  by MIUR ``Programma Dipartimenti di Eccellenza
(2018-2022) - Dipartimento di Matematica F. Casorati
Universit\`a degli Studi di Pavia '' and by INdAM (GNSAGA). In particular, the postdoc position of I. S. was funded by INdAM(GNSAGA).


\begin{thebibliography}{99}
 
	\bibitem{abs} E. Arbarello, A. Bruno, and B. Sernesi. {\em On hyperplane sections of K3 surfaces.}  Algebr. Geom. 4 (2017), no. 5, 562–596. 
	
	\bibitem{bc} E. Ballico and C. Ciliberto. {\em On Gaussian maps for projective varieties.} Geometry of complex projective varieties (Cetraro, 1990), 35–54,
	Sem. Conf., 9, Mediterranean, Rende, 1993. 
	
	\bibitem{bf} C. Barchielli and P. Frediani. {\em On the first Gaussian map for Prym-canonical line bundles.}  Geom. Dedicata 170 (2014), 289–302.
	
	\bibitem{beau} A. Beauville. {\em Vari\'et\'es Kaeleriennes dont la premi\'ere classe de Chern est nulle.} J. Differential Geom. 18 (1983), 755–784. 
	
	\bibitem{bm} A. Beauville and J. Merindol. {\em Sections hyperplanes des surfaces K3.} 	Duke Math. J. 55 (1987), no. 4, 873–878.
	
	\bibitem{bertram coskun} A.Bertram, I. Coskun. {\em 	The birational geometry of the Hilbert scheme of points on surfaces.} Birational geometry, rational curves, and arithmetic, 15–55,
	Simons Symp., Springer, Cham, (2013).  
	
	\bibitem{ccm} A. Calabri, C. Ciliberto, and R. Miranda. {\em The rank of the second Gaussian map for general
		curves.}  Michigan Math. J. 60 (2011), no. 3, 545–559. 
	
	\bibitem{catanesegot} F.Catanese and L. G\"oettsche. {\em d-very-ample line bundles and embeddings of Hilbert schemes of 0-cycles}.  Manuscripta Math 68 (1990), 337–341.
	
	
	\bibitem{chm} C. Ciliberto, J. Harris, and R. Miranda. {\em On the surjectivity of the Wahl map.} Duke Math. Jour. 57 (1988), 829–858.
	
	
	\bibitem{knutsenetal} C.Ciliberto, T. Dedieu, C. Galati, and A. L. Knutsen. {\em Moduli of curves on Enriques surfaces}.  Adv. Math. 365 (2020), 107010-107052.
	
	\bibitem{cv} C. Ciliberto and A. Verra. {\em On the surjectivity of the Gaussian map for Prym-canonical line bundles on a general curve.}  Geometry of complex projective varieties (Cetraro, 1990), 117–141,
	Sem. Conf., 9, Mediterranean, Rende, 1993. 
	
	
	\bibitem{cf2} E. Colombo and P. Frediani.{\em On the second Gaussian map for curves on a K3 surface.} Nagoya
	Math. J. 199 (2010), 123–136.
	
	
	\bibitem{cf3} E. Colombo and P. Frediani.{\em Prym and second Gaussian map for Prym-canonical line bundles.} Advances in Math. 239 (2013), 47-71.
	
	\bibitem{cf1} E. Colombo and P. Frediani. {\em Some results on the second Gaussian map for curves}. Michigan
	Math. J. 58 (2009), 745–758.
	
	
	
	\bibitem{cfp} E.Colombo, P. Frediani and G. Pareschi. {\em Hyperplane sections of abelian surfaces}. J. Algebraic Geom. 21 (2012), 183-200.
	
	
	\bibitem{fog1} J. Fogarty. {\em Algebraic Families on an Algebraic Surface.}  Amer. J. Math. 90 (1968), 511–521. 
	
	
	\bibitem{fog2} J. Fogarty.  {\em Algebraic Families on an Algebraic Surface, II, the Picard Scheme of the Punctual Hilbert Scheme.} Amer. J. Math. 95 (1973), 660–687. 
	
	
	\bibitem{gotsche} L. Göttsche. {\em 
		Hilbert schemes of zero-dimensional subschemes of smooth varieties.} Lecture Notes in Mathematics 1572, Springer-Verlag, Berlin, (1994). 
	
	
	
	\bibitem{kob} S. Kobayashi. {\em First Chern class and holomorphic tensor fields.} Nagoya Math. J. 77 (1980), 5–11. 
	
	\bibitem{knuno} A. L. Knutsen. {\em On kth-order embeddings of K3 surfaces and Enriques surfaces.} Manuscripta Math. 104 (2001), 211–237.
	
	
	\bibitem{knudue} A. L. Knutsen and A. F. Lopez. {\em Surjectivity of Gaussian maps for curves on Enriques surfaces.} Adv. Geom. 7 (2007), no. 2, 215–247. 
	
	
	\bibitem{knutre} A. L. Knutsen, A. F. Lopez and Roberto Muñoz. {\em On the Extendability of Projective Surfaces and a Genus Bound for Enriques-Fano Threefolds.}  J. Differential Geom. 88 (2011), no. 3, 485–518. 
	
	\bibitem{laz2} R. Lazarsfeld. {\em Positivity in Algebraic Geometry I. Classical Setting: Line bundles and Linear Series} Springer, New York (2004).	
	
	\bibitem{mm} S. Mori and S. Mukai. {\em The uniruledness of the moduli space of curves of genus 11.} Algebraic Geometry Proceedings Tokyo, Kyoto, 1982, Springer LNM 1016 (1983), 334–353.
	
	\bibitem{nuer1} H. Nuer. {\em Projectivity and birational geometry of Bridgeland moduli spaces on an Enriques surface.}.  Proc. Lond. Math. Soc. (3) 113 (2016), no. 3, 345–386. 
	
	
	\bibitem{ortiz} A. D. Rios Ortiz. {\em Higher Gaussian Maps on K3 surfaces}. arXiv:2112.09101 
	
	\bibitem{szem} T. Szemberg. {\em On positivity of line bundles on Enriques surfaces}. Trans. Amer. Math. Soc. 353 (2001), no. 12, 4963–4972. 
	
	
	
	\bibitem{wahl intro} J. Wahl. {\em  Introduction to Gaussian maps on an algebraic curve.} Complex projective geometry (Trieste, 1989/Bergen, 1989), 304–323,
	London Math. Soc. Lecture Note Ser., 179, Cambridge Univ. Press, Cambridge, 1992. 
	
	\bibitem{wahl jac} J. Wahl. {\em The Jacobian algebra of a graded Gorenstein singularity.} 	Duke Math. J. 55 (1987), no. 4, 843–871.
	
\end{thebibliography}


\end{document}